\documentclass[12pt]{article}
\usepackage{pifont,indentfirst}
\usepackage{latexsym, amssymb, amsmath, amscd, amsfonts ,texdraw} 

\pagebreak[2]

\newtheorem{thm}{Theorem}[section]

\def\qed{\hfill \rule{4pt}{7pt}}

\makeatletter \@addtoreset{figure}{section} \makeatother

\begin{document}

\begin{center}
{\large \bf Lattice  Polynomials, $12312$-Avoiding Partial Matchings and

Even Trees}
\end{center}

\begin{center}
{\small
William Y. C. Chen$^{1}$, Louis W. Shapiro$^2$ and Susan Y. J. Wu$^3$ \\[5pt]

$^{1}$Center for Combinatorics, LPMC-TJKLC\\
Nankai University, Tianjin 300071, P. R. China

$^2$Department of Mathematics\\
Howard University, Washington, DC 20059, USA

$^3$College of Mathematical Sciences\\
Tianjin Normal University, Tianjin 300387, P. R. China

$^1$chen@nankai.edu.cn, $^2$lshapiro@howard.edu,
$^3$wuyijun@126.com}
\end{center}

\noindent{\bf Abstract.}
The lattice polynomials $L_{i,j}(x)$ are introduced by Hough and Shapiro
as a weighted count  of certain lattice paths from the origin to
the point $(i,j)$. In particular, $L_{2n, n}(x)$ reduces to the
generating function of the numbers $T_{n,k}={1\over n}{n-1+k\choose n-1}{2n-k\choose n+1}$,
which can be viewed as a refinement of the $3$-Catalan numbers $T_n=\frac{1}{2n+1}{3n\choose n}$.
In this paper, we establish a correspondence
between $12312$-avoiding partial matchings and lattice paths, and we show that
the weighted count of such partial matchings with respect to the number of
crossings in a more general sense coincides with the lattice polynomials $L_{i,j}(x)$.
We also introduce a statistic on even  trees, called the $r$-index,
and show that the number of even trees with $2n$ edges and with $r$-index $k$  equal to $T_{n,k}$.



\section{Introduction}

\begin{figure}[h,t]
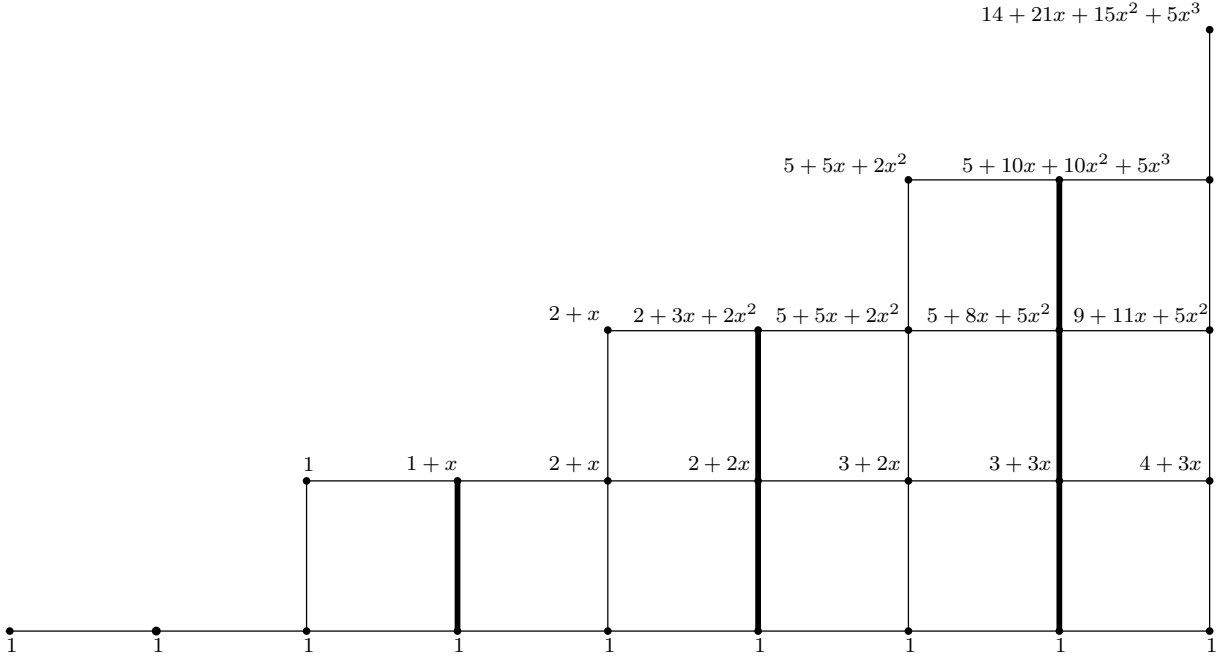

\centertexdraw{
\drawdim mm

\linewd 0.2
\move(0 10) \rlvec(160 0) \rlvec(0 80)
\move(40 10) \rlvec(0 20) \rlvec(120 0)
\move(80 10) \rlvec(0 40) \rlvec(80 0)
\move(120 10) \rlvec(0 60) \rlvec(40 0)
\linewd 0.8
\move(60 10) \rlvec(0 20)
\move(100 10) \rlvec(0 40)
\move(140 10) \rlvec(0 60)

\move(0.53 10) \fcir f:0 r:0.53
\move(20 10) \fcir f:0 r:0.63
\move(40 10) \fcir f:0 r:0.53
\move(60 10) \fcir f:0 r:0.53
\move(80 10) \fcir f:0 r:0.53
\move(100 10) \fcir f:0 r:0.53
\move(120 10) \fcir f:0 r:0.53
\move(140 10) \fcir f:0 r:0.53
\move(160 10) \fcir f:0 r:0.53

\move(40 30) \fcir f:0 r:0.53
\move(60 30) \fcir f:0 r:0.53
\move(80 30) \fcir f:0 r:0.53
\move(100 30) \fcir f:0 r:0.53
\move(120 30) \fcir f:0 r:0.53
\move(140 30) \fcir f:0 r:0.53
\move(160 30) \fcir f:0 r:0.53

\move(80 50) \fcir f:0 r:0.53
\move(100 50) \fcir f:0 r:0.53
\move(120 50) \fcir f:0 r:0.53
\move(140 50) \fcir f:0 r:0.53
\move(160 50) \fcir f:0 r:0.53

\move(120 70) \fcir f:0 r:0.53
\move(140 70) \fcir f:0 r:0.53
\move(160 70) \fcir f:0 r:0.53

\move(160 90) \fcir f:0 r:0.53

\drawdim mm \linewd 0.1

\move(1.5 8) \textref h:R v:C \htext{\scriptsize{$1$}}
\move(21 8) \textref h:R v:C \htext{\scriptsize{$1$}}
\move(41 8) \textref h:R v:C \htext{\scriptsize{$1$}}
\move(61 8) \textref h:R v:C \htext{\scriptsize{$1$}}
\move(81 8) \textref h:R v:C \htext{\scriptsize{$1$}}
\move(101 8) \textref h:R v:C \htext{\scriptsize{$1$}}
\move(121 8) \textref h:R v:C \htext{\scriptsize{$1$}}
\move(141 8) \textref h:R v:C \htext{\scriptsize{$1$}}
\move(161 8) \textref h:R v:C \htext{\scriptsize{$1$}}

\move(41 32) \textref h:R v:C \htext{\scriptsize{$1$}}
\move(60 32) \textref h:R v:C \htext{\scriptsize{$1+x$}}
\move(79 32) \textref h:R v:C \htext{\scriptsize{$2+x$}}
\move(99 32) \textref h:R v:C \htext{\scriptsize{$2+2x$}}
\move(119 32) \textref h:R v:C \htext{\scriptsize{$3+2x$}}
\move(139 32) \textref h:R v:C \htext{\scriptsize{$3+3x$}}
\move(159 32) \textref h:R v:C \htext{\scriptsize{$4+3x$}}

\move(79 52) \textref h:R v:C \htext{\scriptsize{$2+x$}}
\move(100 52) \textref h:R v:C \htext{\scriptsize{$2+3x+2x^2$}}
\move(119 52) \textref h:R v:C \htext{\scriptsize{$5+5x+2x^2$}}
\move(139 52) \textref h:R v:C \htext{\scriptsize{$5+8x+5x^2$}}
\move(160 52) \textref h:R v:C \htext{\scriptsize{$9+11x+5x^2$}}

\move(120 72) \textref h:R v:C \htext{\scriptsize{$5+5x+2x^2$}}
\move(155 72) \textref h:R v:C \htext{\scriptsize{$5+10x+10x^2+5x^3$}}

\move(159 92) \textref h:R v:C \htext{\scriptsize{$14+21x+15x^2+5x^3$}}

}\caption{Lattice polynomials $L_{i,j}(x)$.}\label{f2}
\end{figure}

The lattice polynomials $L_{i,j}(x)$
are introduced by Hough and Shapiro \cite{HoSh} as a weighted count of  lattice paths from the origin $(0,0)$
to $(i,j)$ consisting of unit east steps $(1,0)$ and
north steps $(0,1)$ such that no step goes above the line $x=2y$.
To be more specific, a
north step from $(k,\ell)$ to $(k,\ell+1)$ is given a weight $x$ if
$k$ is odd, and the other steps are assumed to have weight 1. The
weight of a path is the product of the weights of all the steps.
In particular, $L_{2n, n}(x)$ reduces to the
generating function of the numbers
\[ T_{n,k}={1\over n}{n-1+k\choose n-1}{2n-k\choose n+1},\]
which can be viewed as a refinement of the $3$-Catalan numbers \[ T_n=\frac{1}{2n+1}{3n\choose n}.\]
 Figure \ref{f2} gives the first values of $L_{i,j}(x)$.
Note that  $L_{i,0}(x)=1$ for any $i\geq 0$.

The lattice polynomials $L_{i,j}(x)$ are related to
 the descent polynomials on noncrossing trees. Let
$d_n(k,j)$ be the number of noncrossing trees on $\{1,2,\ldots,n\}$
with root degree $j$ and  descent number $k$. Then the descent polynomial
is defined as
\begin{equation}
D_n(x,y)= \sum_{k,j}  d_n(k,j) x^k y^j.
\end{equation}
Hough and Shapiro \cite{HoSh} derived the following relation
\begin{align}
L_{2n,n}(x)&=D_{n+1}(x,1),\label{cw1}\\[5pt]
L_{2n-1,j}(x)& = [y^{n-j}]D_{n+1}(x,y),  \ \ \  \mbox{for $0\leq j\leq n-1$},
\end{align}
where $[y^{n-j}]D_{n+1}(x,y)$ stands for the coefficient of
$y^{n-j}$ in $D_{n+1}(x,y)$. Hough \cite{Hough}
showed that
\begin{equation}\label{e1}
D_n(x,1)=\sum_{k=0}^{n-2}{1\over n-1}{n-2+k\choose
n-2}{2n-2-k\choose n}x^k.
\end{equation}
In view of (\ref{cw1}), the above expression for $D_{n+1}(x,1)$ can be
considered as a formula for $L_{2n,n}(x)$.


This paper is motivated by the observation that the polynomials
$D_{n}(x,1)$ also arise in the context of pattern avoiding matchings.
Let $T_{n,k}$ be the number of  $12312$-avoiding matchings on
$\{1,2,\ldots, 2n\}$ with $k$ crossings. Chen, Mansour and Yan
\cite{ChMY} have shown that
\begin{equation} \label{tnk-1}
T_{n,k}={1\over n}{n-1+k\choose n-1}{2n-k\choose n+1}.
\end{equation}Denote by $T_n(x)$  the generating function of $T_{n,k}$, namely,
\[ T_n(x)= \sum_{k=0}^{n-1} T_{n,k}x^k.\]
It is easy to check that
\[D_{n+1}(x,1)=T_n(x),\]
which can be rewritten as
\begin{equation}\label{cw2}
L_{2n,n}(x)=T_n(x).
\end{equation}

 We introduce a class of $12312$-avoiding partial
 matchings with $i$ points,  denoted by $\mathcal{Q}_i(12312)$,
  and we  establish a correspondence between
 lattice paths counted  by the lattice polynomial $L_{i,j}(x)$
 and the partial matchings in $\mathcal{Q}_i(12312)$ counted by
  a polynomial $Q_{i,j}(x)$.
 Let $q_{i,j,k}$ be the number of $12312$-avoiding
 partial matchings in $\mathcal{Q}_i(12312)$ with $j$ edges and $k$ crossings, where
    the number of crossings of a partial matching is defined
    in a more general sense that an isolated point covered
    by an edge is also considered as a crossing.
    Then $Q_{i,j}(x)$ is the generating function of the numbers $q_{i,j,k}$
    summed over $k$.
In particular,  we see that
$Q_{2n,n}(x)$ equals $T_n(x)$, thus our correspondence
 gives a combinatorial interpretation of (\ref{cw2}).

The second result of this paper is concerned with a refined
enumeration of even trees.
The number of even trees with $2n$ edges is known to be the $3$-Catalan
number $T_n$. We introduce the $r$-index of an even tree.
We construct a correspondence between
lattice paths from $(0,0)$ to $(2n,n)$ with $k$ north steps at odd
positions and even trees with $2n$ edges and $r$-index $k$. Thus
the lattice polynomial $L_{2n,n}(x)$ serves as a formula for
the weighted count of even trees with respect to the $r$-index.

\section{$12312$-avoiding partial matchings}

In this section, we introduce the polynomials $Q_{i,j}(x)$
in connection with the enumeration of partial matchings with
$i$ points, $j$ edges and $k$ crossings.
By constructing a bijection, we show that the polynomials $Q_{i,j}(x)$ coincide
with the lattice polynomials $L_{i,j}(x)$.

We first recall some definitions.  A (complete) matching on a set
$[2n]=\{1,2,\ldots,2n\}$ is a graph with vertex set $[2n]$ and
with $n$ disjoint edges. We assume that the vertices or the points
of a matching is arranged on a horizontal line in increasing
order. Denote an edge of a matching by
$e=(i,j)$ with $i<j$. Two edges $e=(i,j)$ and $e'=(i',j')$ form a
crossing if $i<i'<j<j'$. Denote by
$cr(M)$  the number of crossings  of
$M$.

\begin{figure}[h,t]
\begin{picture}(150,30)
\centertexdraw{ \drawdim mm

\linewd 0.2
\move(10 10) \clvec (13 18)(43 18)(46 10)
\move(16 10) \clvec (18 14)(26 14)(28 10)
\move(22 10) \clvec (24 17)(50 17)(52 10)
\move(34 10) \clvec (35 13)(39 13)(40 10)

\move(11 7) \textref h:R v:C \htext{$1$}
\move(17 7) \textref h:R v:C \htext{$2$}
\move(23 7) \textref h:R v:C \htext{$3$}
\move(29 7) \textref h:R v:C \htext{$4$}
\move(35 7) \textref h:R v:C \htext{$5$}
\move(41 7) \textref h:R v:C \htext{$6$}
\move(47 7) \textref h:R v:C \htext{$7$}
\move(53 7) \textref h:R v:C \htext{$8$}

\move(65 12) \textref h:R v:C \htext{$\Longrightarrow$}

\move(70 10) \clvec (73 18)(103 18)(106 10)
\move(76 10) \clvec (78 14)(86 14)(88 10)
\move(82 10) \clvec (84 17)(110 17)(112 10)
\move(94 10) \clvec (95 13)(99 13)(100 10)

\move(71 7) \textref h:R v:C \htext{$1$}
\move(77 7) \textref h:R v:C \htext{$2$}
\move(83 7) \textref h:R v:C \htext{$3$}
\move(89 7) \textref h:R v:C \htext{$2$}
\move(95 7) \textref h:R v:C \htext{$4$}
\move(101 7) \textref h:R v:C \htext{$4$}
\move(107 7) \textref h:R v:C \htext{$1$}
\move(113 7) \textref h:R v:C \htext{$3$}

\move(125 12) \textref h:R v:C \htext{$\Longrightarrow$}

\move(155 12) \textref h:R v:C \htext{$P=12324413$}
}
\end{picture}
\caption{A matching and its canonical sequential form.}\label{f1}
\end{figure}
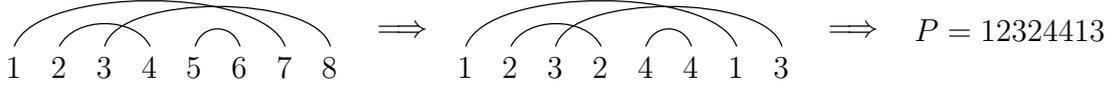

A matching  can be expressed by its canonical sequential form
\cite{Klaz}, or the Davenport-Schinzel sequence \cite{MuSt}.
For a matching $P$ on $[2n]$, denote the edges of $P$ by
$e_1=(i_1,j_1)$, $e_2=(i_2,j_2)$, \ldots, $e_n=(i_n,j_n)$, where
$i_1<i_2<\cdots<i_n$. We write $P=a_1a_2\cdots a_{2n}$, where $a_i=j$
if $i$ is an endpoint of the edge $e_j$.
 For example, the matching in Figure \ref{f1}
can be expressed by the sequence $12324413$.

Let $\alpha=\alpha_1\alpha_2\cdots\alpha_k$ and
$\pi=\pi_1\pi_2\cdots\pi_k$ be two sequences. We say $\alpha$ and
$\pi$ are order-isomorphic if for any $1\leq i,j\leq k$,
$\alpha_i<\alpha_j$ (resp. $\alpha_i=\alpha_j$, $\alpha_i>\alpha_j$)  if and only if $\pi_i<\pi_j$ (resp. $\pi_i=\pi_j$, $\pi_i>\pi_j$).   We say that a
canonical sequential form $P$  avoids a sequence
 $\pi$, or $P$ is $\pi$-avoiding,
if no subsequence of $P$ is order-isomorphic to $\pi$. Such a sequence
$\pi$
is usually called a  pattern. Denote by $\mathcal{M}_n(\pi)$ the set
of  matchings on $[2n]$ which avoid a pattern $\pi$. It has
been shown by Chen, Mansour and Yan \cite{ChMY} that the number of
$12312$-avoiding matchings on $[2n]$ equals the $3$-Catalan number, namely,
$|\mathcal{M}_n(12312)|=T_n.$

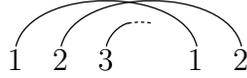
\begin{figure}[h,t]
\begin{picture}(150,20)
\centertexdraw{ \drawdim mm

\linewd 0.2
\move(10 10) \clvec (12 18)(32 18)(34 10)
\move(16 10) \clvec (18 18)(38 18)(40 10)
\move(22 10) \clvec (22.5 12)(24.5 13)(25 13)
\lpatt(0.3 0.5)
\move(25 13) \clvec (26 13.2)(27 13.2)(28 13)

\move(11 8) \textref h:R v:C \htext{$1$}
\move(17 8) \textref h:R v:C \htext{$2$}
\move(23 8) \textref h:R v:C \htext{$3$}
\move(35 8) \textref h:R v:C \htext{$1$}
\move(41 8) \textref h:R v:C \htext{$2$}
}
\end{picture}
\caption{The  pattern 12312. }\label{bad-pattern}
\end{figure}

\begin{figure}[h,t]
\begin{picture}(150,25)
\centertexdraw{ \drawdim mm

\linewd 0.2
\move(10 10) \clvec (12 18)(28 18)(28 10)
\move(16 10) \clvec (18 18)(32 18)(34 10)
\move(22 10) \fcir f:0 r:0.3

\move(50 10) \clvec (52 18)(68 18)(68 10)
\move(62 10) \clvec (64 17)(72 17)(74 10)
\move(56 10) \fcir f:0 r:0.3

\move(90 10) \clvec (92 18)(112 18)(114 10)
\move(102 10) \clvec (104 15)(106 15)(108 10)
\move(96 10) \fcir f:0 r:0.3

\move(11 7) \textref h:R v:C \htext{\footnotesize{$1$}}
\move(17 7) \textref h:R v:C \htext{\footnotesize{$2$}}
\move(23 7) \textref h:R v:C \htext{\footnotesize{$3$}}
\move(29 7) \textref h:R v:C \htext{\footnotesize{$1$}}
\move(35 7) \textref h:R v:C \htext{\footnotesize{$2$}}

\move(51 7) \textref h:R v:C \htext{\footnotesize{$1$}}
\move(57 7) \textref h:R v:C \htext{\footnotesize{$2$}}
\move(63 7) \textref h:R v:C \htext{\footnotesize{$3$}}
\move(69 7) \textref h:R v:C \htext{\footnotesize{$1$}}
\move(75 7) \textref h:R v:C \htext{\footnotesize{$3$}}

\move(91 7) \textref h:R v:C \htext{\footnotesize{$1$}}
\move(97 7) \textref h:R v:C \htext{\footnotesize{$2$}}
\move(103 7) \textref h:R v:C \htext{\footnotesize{$3$}}
\move(109 7) \textref h:R v:C \htext{\footnotesize{$3$}}
\move(115 7) \textref h:R v:C \htext{\footnotesize{$1$}}

}
\end{picture}
\caption{The partial matchings. }\label{partial-pattern}
\end{figure}
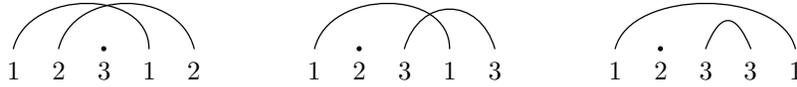

A partial matching on $[m]=\{1, 2,\ldots, m\}$ can be viewed as a
matching on a subset of $[m]$. For example,  Figure \ref{partial-pattern}
 gives three partial matchings with
$5$ vertices.

We shall consider a special class of partial matchings such that any
edge covers at most one isolated point and the number of
isolated points to the left of any edge cannot exceed  the
number of isolated points to the right of this edge. Denote this
class of partial matchings on $[m]$ by $\mathcal{Q}_m$. For a
partial matching $M\in\mathcal{Q}_m$, an edge $e=(i,j)$ and an
isolated point $k$ form a crossing if $i<k<j$. Denoted by
$cr(M)$ the number of crossings of $M$, where we count the
crossings formed by two edges, as well as one edge and an
isolated point. For example, for the partial
matchings in Figure \ref{partial-pattern}, we have
 $cr(12312)=3$,
$cr(12313)=2$, and $cr(12331)=1$.

Let $\mathcal{Q}_m(\pi)$ denote the set of
 partial matchings in $\mathcal{Q}_m$ that avoid pattern
$\pi$. We shall be concerned with a refined
enumeration of the set $\mathcal{Q}_i(12312)$.
 Let  $q_{i,j,k}$ be the number of partial matchings in the set $\mathcal{Q}_{i}(12312)$
with $j$ edges and $k$ crossings. Set
\begin{equation}
Q_{i,j}(x) = \sum_k  q_{i,j,k}\,x^k,
\end{equation}
where $k$ ranges from $0$ to $\lceil{i-1\over 2}\rceil$.
In this notation, we have the following relation.

\begin{thm} The  polynomial $Q_{i,j}(x)$ equals
the lattice polynomial $L_{i,j}(x)$.
\end{thm}

\noindent{\it Proof.} We proceed to give a procedure to generate all
partial matchings in $\mathcal{Q}_m(12312)$.
For this purpose, we introduce two operations, called the shifting
operation and the lifting operation on partial matchings
corresponding to the east and north steps in the lattice path.

The shifting operation is defined by adding an isolated point to
the right  of the last vertex of $M$. This operation corresponds to an
east step $E=(1,0)$ in the lattice path.

The lifting operation is defined by adding an edge in the middle of
the partial matching $M$. More precisely, if there are $2k$ isolated
points in $M$, then connect the $k$-th and $(k+1)$-st isolated points to form a new edge. If there are $2k+1$ isolated points in $M$, then connect the $k$-th and $(k+2)$-nd isolated points. This operation corresponds to a north step $N=(0,1)$ in
the lattice path.

We now describe the procedure to generate
all  partial matchings in $\mathcal{Q}_m(12312)$.
Start with an empty set $\emptyset$ at the origin
$(0,0)$. Suppose that we have constructed a partial matching
$M\in\mathcal{Q}_m(12312)$ at the position $(i,j)$. Using the
shifting operation or the lifting operation on $M$, we
obtain a new partial matching $M'$ at the position $(i+1,j)$ or $(i,j+1)$.
Iterating this process, we generate all the partial matchings at
a given position $(i',j')$.

\begin{figure}[h,t]
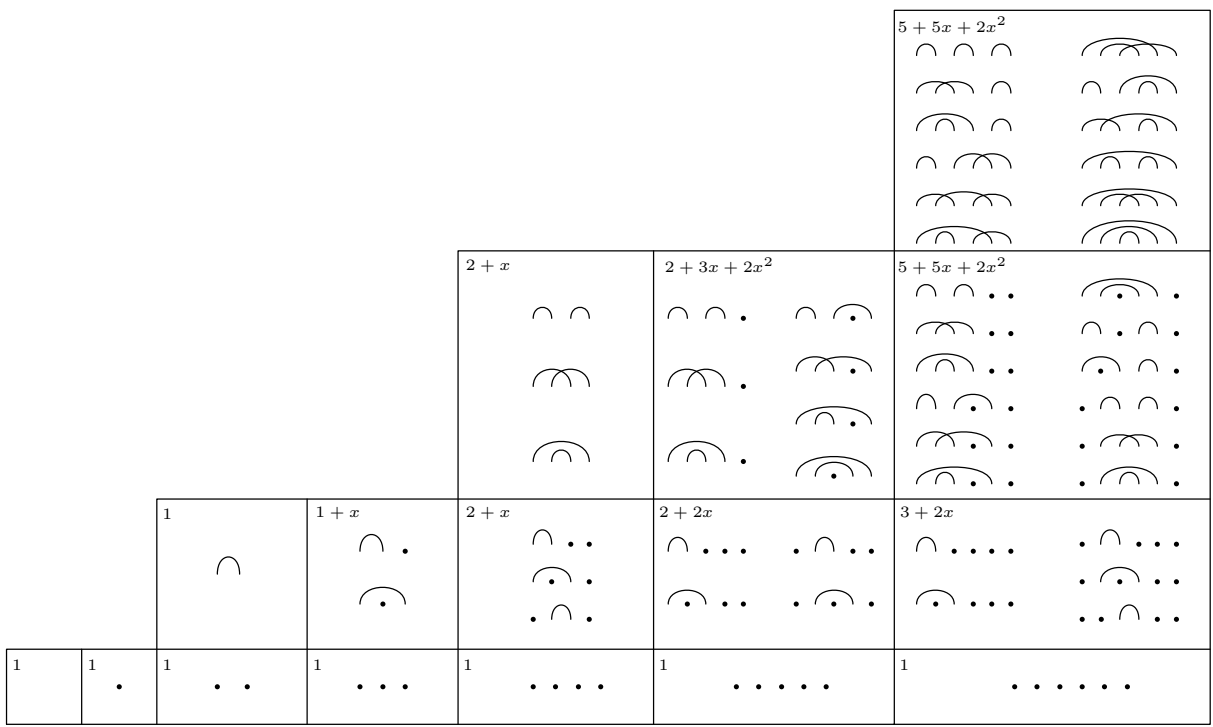

\centertexdraw {

\drawdim mm \linewd 0.2


\setgray 0 \move(0 0) \rlvec(160 0) \rlvec(0 95) \rlvec(-42 0)
\rlvec(0 -95) \move(0 0) \rlvec(0 10) \rlvec(160 0) \move(10 0)
\rlvec(0 10) \move(20 0) \rlvec(0 30) \rlvec(140 0) \move(40 0)
\rlvec(0 30) \move(60 0) \rlvec(0 63) \rlvec(100 0) \move(86 0)
\rlvec(0 63)

\move(15 5) \fcir f:0 r:0.3

\move(28 5) \fcir f:0 r:0.3 \move(32 5) \fcir f:0 r:0.3

\move(47 5) \fcir f:0 r:0.3 \move(50 5) \fcir f:0 r:0.3 \move(53 5)
\fcir f:0 r:0.3

\move(70 5) \fcir f:0 r:0.3 \move(73 5) \fcir f:0 r:0.3 \move(76 5)
\fcir f:0 r:0.3 \move(79 5) \fcir f:0 r:0.3

\move(97 5) \fcir f:0 r:0.3 \move(100 5) \fcir f:0 r:0.3 \move(103
5) \fcir f:0 r:0.3 \move(106 5) \fcir f:0 r:0.3 \move(109 5) \fcir
f:0 r:0.3

\move(134 5) \fcir f:0 r:0.3 \move(137 5) \fcir f:0 r:0.3 \move(140
5) \fcir f:0 r:0.3 \move(143 5) \fcir f:0 r:0.3 \move(146 5) \fcir
f:0 r:0.3 \move(149 5) \fcir f:0 r:0.3

\linewd 0.15

\move(28 20) \clvec (28 23)(31 23)(31 20)

\move(47 23) \clvec (47 26)(50 26)(50 23)  \move(53 23) \fcir f:0
r:0.3 \move(47 16) \clvec (47 19)(53 19)(53 16)  \move(50 16) \fcir
f:0 r:0.3

\move(70 24) \clvec (70 26.5)(72.5 26.5)(72.5 24) \move(75 24) \fcir
f:0 r:0.3 \move(77.5 24) \fcir f:0 r:0.3 \move(70 19) \clvec (70
21.5)(75 21.5)(75 19) \move(72.5 19) \fcir f:0 r:0.3 \move(77.5 19)
\fcir f:0 r:0.3 \move(72.5 14) \clvec (72.5 16.5)(75 16.5)(75 14)
\move(70 14) \fcir f:0 r:0.3 \move(77.5 14) \fcir f:0 r:0.3

\move(70 54) \clvec (70 56)(72.5 56)(72.5 54) \move(75 54) \clvec
(75 56)(77.5 56)(77.5 54) \move(70 45) \clvec (70 48)(75 48)(75 45)
\move(72.5 45) \clvec (72.5 48)(77.5 48)(77.5 45) \move(72.5 35)
\clvec (72.5 37)(75 37)(75 35) \move(70 35) \clvec (70 38.5)(77.5
38.5)(77.5 35)

\move(88 54) \clvec (88 56)(90.5 56)(90.5 54) \move(93 54) \clvec
(93 56)(95.5 56)(95.5 54) \move(98 54) \fcir f:0 r:0.3

\move(88 45) \clvec (88 48)(93 48)(93 45)  \move(90.5 45) \clvec
(90.5 48)(95.5 48)(95.5 45) \move(98 45) \fcir f:0 r:0.3

\move(90.5 35) \clvec (90.5 37)(93 37)(93 35) \move(88 35) \clvec
(88 38.5)(95.5 38.5)(95.5 35) \move(98 35) \fcir f:0 r:0.3

\move(105 54) \clvec (105 56)(107.5 56)(107.5 54) \move(110 54)
\clvec (110 56.5)(115 56.5)(115 54) \move(112.5 54) \fcir f:0 r:0.3

\move(105 47) \clvec (105 49.5)(110 49.5)(110 47)  \move(107.5 47)
\clvec (107.5 49.5)(115 49.5)(115 47) \move(112.5 47) \fcir f:0
r:0.3

\move(105 40) \clvec (105 43)(115 43)(115 40)  \move(107.5 40)
\clvec (107.5 42)(110 42)(110 40) \move(112.5 40) \fcir f:0 r:0.3

\move(105 33) \clvec (105 36.5)(115 36.5)(115 33) \move(107.5 33)
\clvec (107.5 35.5)(112.5 35.5)(112.5 33) \move(110 33) \fcir f:0
r:0.3

\move(88 23) \clvec (88 25.5)(90.5 25.5)(90.5 23) \move(93 23) \fcir
f:0 r:0.3 \move(95.5 23) \fcir f:0 r:0.3 \move(98 23) \fcir f:0
r:0.3

\move(88 16) \clvec (88 18.5)(93 18.5)(93 16) \move(90.5 16) \fcir
f:0 r:0.3 \move(95.5 16) \fcir f:0 r:0.3 \move(98 16) \fcir f:0
r:0.3

\move(107.5 23) \clvec (107.5 25.5)(110 25.5)(110 23) \move(105 23)
\fcir f:0 r:0.3 \move(112.5 23) \fcir f:0 r:0.3 \move(115 23) \fcir
f:0 r:0.3

\move(107.5 16) \clvec (107.5 18.5)(112.5 18.5)(112.5 16) \move(105
16) \fcir f:0 r:0.3 \move(110 16) \fcir f:0 r:0.3 \move(115 16)
\fcir f:0 r:0.3

\move(121 23) \clvec (121 25.5)(123.5 25.5)(123.5 23) \move(126 23)
\fcir f:0 r:0.3 \move(128.5 23) \fcir f:0 r:0.3 \move(131 23) \fcir
f:0 r:0.3 \move(133.5 23) \fcir f:0 r:0.3

\move(121 16) \clvec (121 18.5)(126 18.5)(126 16) \move(123.5 16)
\fcir f:0 r:0.3 \move(128.5 16) \fcir f:0 r:0.3 \move(131 16) \fcir
f:0 r:0.3 \move(133.5 16) \fcir f:0 r:0.3

\move(145.5 24) \clvec (145.5 26.5)(148 26.5)(148 24) \move(143 24)
\fcir f:0 r:0.3 \move(150.5 24) \fcir f:0 r:0.3 \move(153 24) \fcir
f:0 r:0.3 \move(155.5 24) \fcir f:0 r:0.3

\move(145.5 19) \clvec (145.5 21.5)(150.5 21.5)(150.5 19) \move(143
19) \fcir f:0 r:0.3 \move(148 19) \fcir f:0 r:0.3 \move(153 19)
\fcir f:0 r:0.3 \move(155.5 19) \fcir f:0 r:0.3

\move(148 14) \clvec (148 16.5)(150.5 16.5)(150.5 14) \move(143 14)
\fcir f:0 r:0.3 \move(145.55 14) \fcir f:0 r:0.3 \move(153 14) \fcir
f:0 r:0.3 \move(155.5 14) \fcir f:0 r:0.3

\move(121 57) \clvec (121 59)(123.5 59)(123.5 57) \move(126 57)
\clvec (126 59)(128.5 59)(128.5 57) \move(131 57) \fcir f:0 r:0.3
\move(133.5 57) \fcir f:0 r:0.3

\move(121 52) \clvec (121 54)(126 54)(126 52) \move(123.5 52) \clvec
(123.5 54)(128.5 54)(128.5 52) \move(131 52) \fcir f:0 r:0.3
\move(133.5 52) \fcir f:0 r:0.3

\move(121 47) \clvec (121 50)(128.5 50)(128.5 47) \move(123.5 47)
\clvec (123.5 49)(126 49)(126 47) \move(131 47) \fcir f:0 r:0.3
\move(133.5 47) \fcir f:0 r:0.3

\move(121 42) \clvec (121 44.5)(123.5 44.5)(123.5 42) \move(126 42)
\clvec (126 44.5)(131 44.5)(131 42) \move(128.5 42) \fcir f:0 r:0.3
\move(133.5 42) \fcir f:0 r:0.3

\move(121 37) \clvec (121 39.5)(126 39.5)(126 37) \move(123.5 37)
\clvec (123.5 39.5)(131 39.5)(131 37) \move(128.5 37) \fcir f:0
r:0.3 \move(133.5 37) \fcir f:0 r:0.3

\move(121 32) \clvec (121 35)(131 35)(131 32) \move(123.5 32) \clvec
(123.5 34)(126 34)(126 32) \move(128.5 32) \fcir f:0 r:0.3
\move(133.5 32) \fcir f:0 r:0.3
\move(143 57) \clvec (143 60)(153 60)(153 57) \move(145.5 57) \clvec
(145.5 59)(150.5 59)(150.5 57) \move(148 57) \fcir f:0 r:0.3
\move(155.5 57) \fcir f:0 r:0.3

\move(143 52) \clvec (143 54)(145.5 54)(145.5 52) \move(150.5 52)
\clvec (150.5 54)(153 54)(153 52) \move(148 52) \fcir f:0 r:0.3
\move(155.5 52) \fcir f:0 r:0.3

\move(143 47) \clvec (143 49.5)(148 49.5)(148 47) \move(150.5 47)
\clvec (150.5 49)(153 49)(153 47) \move(145.5 47) \fcir f:0 r:0.3
\move(155.5 47) \fcir f:0 r:0.3

\move(145.5 42) \clvec (145.5 44)(148 44)(148 42) \move(150.5 42)
\clvec (150.5 44)(153 44)(153 42) \move(143 42) \fcir f:0 r:0.3
\move(155.5 42) \fcir f:0 r:0.3

\move(145.5 37) \clvec (145.5 39)(150.5 39)(150.5 37) \move(148 37)
\clvec (148 39)(153 39)(153 37) \move(143 37) \fcir f:0 r:0.3
\move(155.5 37) \fcir f:0 r:0.3

\move(145.5 32) \clvec (145.5 35)(153 35)(153 32) \move(148 32)
\clvec (148 34)(150.5 34)(150.5 32) \move(143 32) \fcir f:0 r:0.3
\move(155.5 32) \fcir f:0 r:0.3

\move(121 89) \clvec (121 91)(123.5 91)(123.5 89) \move(126 89)
\clvec (126 91)(128.5 91)(128.5 89) \move(131 89) \clvec (131
91)(133.5 91)(133.5 89)

\move(121 84) \clvec (121 86)(126 86)(126 84) \move(123.5 84) \clvec
(123.5 86)(128.5 86)(128.5 84) \move(131 84) \clvec (131 86)(133.5
86)(133.5 84)

\move(121 79) \clvec (121 82)(128.5 82)(128.5 79) \move(123.5 79)
\clvec (123.5 81)(126 81)(126 79) \move(131 79) \clvec (131
81)(133.5 81)(133.5 79)

\move(121 74) \clvec (121 76)(123.5 76)(123.5 74) \move(126 74)
\clvec (126 76.5)(131 76.5)(131 74) \move(128.5 74) \clvec (128.5
76.5)(133.5 76.5)(133.5 74)

\move(121 69) \clvec (121 71)(126 71)(126 69) \move(123.5 69) \clvec
(123.5 71.5)(131 71.5)(131 69) \move(128.5 69) \clvec (128.5
71)(133.5 71)(133.5 69)

\move(121 64) \clvec (121 67)(131 67)(131 64) \move(123.5 64) \clvec
(123.5 66)(126 66)(126 64) \move(128.5 64) \clvec (128.5 66)(133.5
66)(133.5 64)

\move(143 89) \clvec (143 92)(153 92)(153 89) \move(145.5 89) \clvec
(145.5 91)(150.5 91)(150.5 89) \move(148 89) \clvec (148 91)(155.5
91)(155.5 89)

\move(143 84) \clvec (143 86)(145.5 86)(145.5 84) \move(150.5 84)
\clvec (150.5 86)(153 86)(153 84) \move(148 84) \clvec (148
87)(155.5 87)(155.5 84)

\move(143 79) \clvec (143 81)(148 81)(148 79) \move(150.5 79) \clvec
(150.5 81)(153 81)(153 79) \move(145.5 79) \clvec (145.5 82)(155.5
82)(155.5 79)

\move(145.5 74) \clvec (145.5 76)(148 76)(148 74) \move(150.5 74)
\clvec (150.5 76)(153 76)(153 74) \move(143 74) \clvec (143
77)(155.5 77)(155.5 74)

\move(145.5 69) \clvec (145.5 71)(150.5 71)(150.5 69) \move(148 69)
\clvec (148 71)(153 71)(153 69) \move(143 69) \clvec (143 72)(155.5
72)(155.5 69)

\move(145.5 64) \clvec (145.5 67)(153 67)(153 64) \move(148 64)
\clvec (148 66)(150.5 66)(150.5 64) \move(143 64) \clvec (143
68)(155.5 68)(155.5 64)

\move(2 8) \textref h:R v:C \htext{\tiny{$1$}} \move(12 8) \textref
h:R v:C \htext{\tiny{$1$}} \move(22 8) \textref h:R v:C
\htext{\tiny{$1$}} \move(42 8) \textref h:R v:C \htext{\tiny{$1$}}
\move(62 8) \textref h:R v:C \htext{\tiny{$1$}} \move(88 8) \textref
h:R v:C \htext{\tiny{$1$}} \move(120 8) \textref h:R v:C
\htext{\tiny{$1$}}

\move(22 28) \textref h:R v:C \htext{\tiny{$1$}} \move(47 28)
\textref h:R v:C \htext{\tiny{$1+x$}} \move(67 28) \textref h:R v:C
\htext{\tiny{$2+x$}} \move(94 28) \textref h:R v:C
\htext{\tiny{$2+2x$}} \move(126 28) \textref h:R v:C
\htext{\tiny{$3+2x$}}

\move(67 61) \textref h:R v:C \htext{\tiny{$2+x$}} \move(102 61)
\textref h:R v:C \htext{\tiny{$2+3x+2x^2$}} \move(133 61) \textref
h:R v:C \htext{\tiny{$5+5x+2x^2$}}

\move(133 93) \textref h:R v:C \htext{\tiny{$5+5x+2x^2$}}

} \caption{ $12312$-avoiding partial matchings and lattice polynomials.}\label{f3}
\end{figure}

Figure \ref{f3} gives an illustration of the above procedure to
generate partial matchings.

Clearly, for a shifting operation, adding a new isolated point to the
rightmost position does not change the number of crossings, i.e.,
$cr(M')=cr(M)$. For a lifting operation, if there are $2k$ isolated
points in $M$, it is easily checked that
connecting the $k$-th and $(k+1)$-st isolated points does not create any crossings, that is,
$cr(M')=cr(M)$. If there are $2k+1$
isolated points in $M$, then connecting the $k$-th and $(k+2)$-nd
isolated points creates a new crossing on $M'$, that is,
$cr(M')=cr(M)+1$.

To show that $Q_{i,j}(x)$ equals the lattice polynomial
$L_{i,j}(x)$, we need to verify the following facts.

Claim 1: After the shifting or lifting operation, the new partial
matching $M'$ still avoids the pattern $12312$.

Claim 2: The two operations do not lead to the same partial matching.

Since the shifting operation is defined by just adding a new isolated
point, it is obvious that $M'$ still avoids the pattern $12312$ after
this operation.

Let us  consider the  lifting operation.
Suppose to the contrary that
$M'$ is not $12312$-avoiding, that is,  there exists  a
subsequence of $M'$ which is of pattern  $12312$.
Assume that the new edge produced by lifting operation in $M'$ is $e$.
 Since $M$ is $12312$-avoiding, the pattern
$12312$ in $M'$ must involve the edge $e$. We have  two cases.

Case 1. There are $2k$ isolated points in $M$.
For any subsequence of $M'$ which is of pattern
$12312$, the maximum point  cannot be an isolated
point; Otherwise there exits an edge in $M$  that covers more
than one point, contradicting  the definition of $\mathcal{Q}_m$.
Hence the subsequence of $M'$ involves three edges $(i_1,j_1)$, $(i_2,j_2)$, $(i_3,j_3)$, where $i_1<i_2<i_3<j_1<j_2$.

If $e=(i_1,j_1)$, according to the position of $j_3$, there are three possibilities as shown in the following figure. Obviously, the number of isolated points to the left of  $(i_2,j_2)$ is greater than the
number of isolated points to the right of the edge $(i_2, j_2)$, which contradicts the definition of $\mathcal{Q}_m$.

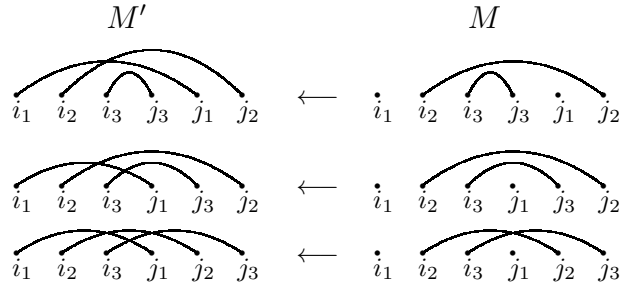
\begin{figure}[h,t]
\setlength{\unitlength}{0.6mm}
\begin{center}
\begin{picture}(150,50)

\put(20,50){$M'$} \put(100,50){$M$}

 \multiput(0,0)(10,0){6}{\circle*{1}}

 \qbezier[1000](0,0)(15,10)(30,0)
\qbezier[1000](10,0)(25,10)(40,0) \qbezier[1000](20,0)(35,10)(50,0)

\multiput(80,0)(10,0){6}{\circle*{1}}

 \qbezier[1000](90,0)(105,10)(120,0)
\qbezier[1000](100,0)(115,10)(130,0)

\footnotesize

\put(-1,-5){$i_1$}\put(9,-5){$i_2$}\put(19,-5){$i_3$}
\put(29,-5){$j_1$}\put(39,-5){$j_2$}\put(49,-5){$j_3$}

\put(79,-5){$i_1$}\put(89,-5){$i_2$}\put(99,-5){$i_3$}
\put(109,-5){$j_1$}\put(119,-5){$j_2$}\put(129,-5){$j_3$}

\put(62,-2){$\longleftarrow$}


\multiput(0,15)(10,0){6}{\circle*{1}}

 \qbezier[1000](0,15)(15,25)(30,15)
\qbezier[1000](10,15)(30,30)(50,15)
\qbezier[1000](20,15)(30,25)(40,15)

\multiput(80,15)(10,0){6}{\circle*{1}}

 \qbezier[1000](90,15)(110,30)(130,15)
\qbezier[1000](100,15)(110,25)(120,15)

\footnotesize

\put(-1,10){$i_1$}\put(9,10){$i_2$}\put(19,10){$i_3$}
\put(29,10){$j_1$}\put(39,10){$j_3$}\put(49,10){$j_2$}

\put(79,10){$i_1$}\put(89,10){$i_2$}\put(99,10){$i_3$}
\put(109,10){$j_1$}\put(119,10){$j_3$}\put(129,10){$j_2$}

\put(62,13){$\longleftarrow$}


\multiput(0,35)(10,0){6}{\circle*{1}}

 \qbezier[1000](0,35)(20,50)(40,35)
\qbezier[1000](10,35)(30,55)(50,35)
\qbezier[1000](20,35)(25,45)(30,35)

\multiput(80,35)(10,0){6}{\circle*{1}}

 \qbezier[1000](90,35)(110,50)(130,35)
\qbezier[1000](100,35)(105,45)(110,35)

\footnotesize

\put(-1,30){$i_1$}\put(9,30){$i_2$}\put(19,30){$i_3$}
\put(29,30){$j_3$}\put(39,30){$j_1$}\put(49,30){$j_2$}

\put(79,30){$i_1$}\put(89,30){$i_2$}\put(99,30){$i_3$}
\put(109,30){$j_3$}\put(119,30){$j_1$}\put(129,30){$j_2$}

\put(62,33){$\longleftarrow$}
\end{picture}
\end{center}
\end{figure}

The same argument applies to the situations
 $e=(i_2,j_2)$ and  $e=(i_3,j_3)$.

Case 2. There are $2k+1$ isolated points in $M$.
It is easily seen that
any subsequence of $M'$ which is of pattern $12312$
consists of two edges and one isolated point.
If $e=(i_1,j_1)$ or $e=(i_2,j_2)$,
there is always one edge of $M$ covering more than one isolated point,
which contradicts the definition of  $\mathcal{Q}_m$. So
 Claim 1 is proved.

We turn to the proof of Claim 2.
We use induction on $i$ and $j$. It is easy to check the
statement  holds for small values of $(i,j)$, see
Figure \ref{f3}. We assume that the claim holds for the positions
$(i-1, j)$ and $(i,j-1)$.

Suppose  that $M_1'$ and $M_2'$ are
$12312$-avoiding partial matchings in the  position $(i,j)$,
and that $M_1'$ is obtained from $M_1$ and $M_2'$ is obtained from $M_2$.
There are three cases for the positions of $M_1$ and $M_2$.

If both $M_1$ and $M_2$ are in the position $(i-1,j)$ or both $M_1$ and $M_2$ are
in the position $(i,j-1)$, then by the induction hypothesis $M_1$ and
$M_2$ are distinct. It follows that  $M_1'$ and $M_2'$
are distinct.

We now consider the case that $M_1$ is in the position $(i-1,j)$ and
$M_2$ is in the position $(i,j-1)$. The assumption implies that
  $i,j\geq 1$. Clearly, $M_1'$ is obtained from $M_1$ via a shifting operation by adding
an isolated point to the rightmost point of $M_1$.
Evidently, the number
of isolated points to the left of any edge in $M_1'$
is less than the number of isolated points to the right of this edge.
Since $M_2'$ is obtained from $M_2$ via a lifting operation by
adding an edge, we see that the number of isolated points to the left of this edge is equal to the number of isolated points to the right of this edge in $M_2'$.
Therefore, in this case $M_1'$ and $M_2'$ are distinct.
 So we arrive at the conclusion that
 the partial matchings obtained at the position $(i,j)$ are distinct.
 Thus the claim holds by induction.

It remains to show that the above procedure generates all
partial matchings in $\mathcal{Q}_i$ with $j$ edges.
To this end, we give the reverse procedure to construct a
lattice path from any partial matching in $\mathcal{Q}_i$ with $j$ edges.

Let $M$ be a $12312$-avoiding partial matching in $\mathcal{Q}_i$ with  $j$ edges.
The construction is recursive.
More precisely, we shall construct a lattice
path from $(i,j)$ to the origin that will not go beyond the line $x = 2y$.
We start at the position  $(i,j)$. If  for any edge in $M$ the number of isolated points to the left of this edge is less than the number of isolated points to the right of this edge, then there must be an isolated point in the end. We delete  the last isolated point and continue to consider the
  position $(i-1,j)$. Note that in this case $i\geq 2j+1$ since there are isolated points in $M$, and thus $(i-1,j)$ does  not go beyond the line $x=2y$. Otherwise, among the edges in $M$ for which the number of isolated points to the left equals the number of isolated points to the right, we can choose the  unique edge $e$ with the rightmost right endpoint. Then, we delete the edge $e$ and continue to
  consider the position $(i,j-1)$. Obviously, $(i,j-1)$ does not exceed the line $x=2y$. Iterating this procedure, we eventually get the required lattice path.
\qed

The proof of the above theorem can be considered as a recursive construction of
a correspondence between lattice paths counted by $L_{i,j}(x)$ and $12312$-avoiding
partial matchings counted by $Q_{i,j}(x)$.

\section{The $r$-index of even trees}

In this section, we define a statistic, which we call the $r$-index,  on
even  trees. We shall show that the generating function for the number
of even trees with $2n$ edges and  with $r$-index equal to $k$
coincides with the lattice polynomial $L_{2n,n}(x)$.

Recall that an even tree is a plane tree in which every vertex has
an even number of children. The number of even trees with $2n$ edges
equals the $3$-Catalan number $T_n$. 
If a vertex has $2k$ children, we call the first $k$ children left
children and the last $k$ children right children. So every vertex except for the
root is either a left child or a right child. 
The {$r$-index} of
an even tree $T$,  denoted
by $r(T)$, is defined as half of the sum of the degrees
 of right children, where the degree of a node is meant to be the
number of its children. Define $R_{n,k}$ to be the number of even trees
with $2n$ edges and with the $r$-index equal to $k$, and define the
generating function of $R_{n,k}$ by
\begin{equation}
R_n(x)=\sum_{k=0}^{n-1} R_{n,k} x^k.
\end{equation}

The following theorem gives a connection between
lattice polynomials and even trees counted with respect to
the number of edges and the $r$-index. 
The proof can be viewed as a recursive construction of
a bijection.

\begin{thm}\label{thm2.1}
For $n\geq 1$,  $R_n(x)$ equals the lattice polynomial $L_{2n,n}(x)$.
\end{thm}

\noindent {\it Proof.} We proceed to present
 a procedure to generate
even trees  parallel to the construction of
 lattice paths counted
by the lattice  polynomials. To this end, we shall introduce two
operations, called the shifting operation and the lifting operation,
on even trees corresponding to the east and north steps of the
lattice paths. These two operations can be viewed as the actions to
shift the even trees in the lattice from left to right, and to lift
the even trees.

Actually, we need an intermediate structure for the generation of
even trees that corresponds to the intermediate points, namely, the
points with even $x$-coordinates,  for the generation of lattice  polynomials. Such an intermediate structure will be represented by
an even tree with a pair of dotted edges from the root to its first
child and to its last child.

We first describe the shifting operation. We start with an empty set
$\emptyset$ at origin  $(0,0)$. Let  $T$ be an even tree
which may contain a pair of dotted edges as the outside edges of the
root. The shifting operation is meant to transform  $T$ in a
position $(i,j)$ to an even tree $T'$ in the position
to the immediately right of $T$. There are two cases.

 If $T$
contains a pair of dotted edges, then $T'$ is obtained from $T$ by
changing the dotted edges to regular edges. If $T$ does not contain
any dotted edges, then $T'$ is obtained from $T$ by adding a pair of
dotted edges as the outside edges of the root.

\begin{figure}[h,t]
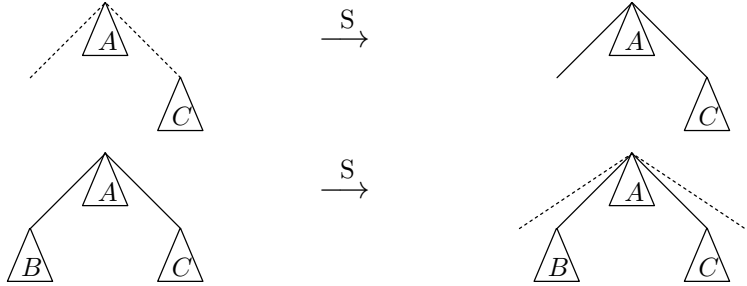

\centertexdraw{
\drawdim mm

\linewd 0.2

\lpatt(0.3 0.5)
\move(10 60) \rlvec(10 10) \rlvec(10 -10)
\lpatt()
\move(80 60) \rlvec(10 10) \rlvec(10 -10)

\move(20 70) \rlvec(3 -7) \rlvec(-6 0) \rlvec(3 7)
\move(90 70) \rlvec(3 -7) \rlvec(-6 0) \rlvec(3 7)
\move(30 60) \rlvec(3 -7) \rlvec(-6 0) \rlvec(3 7)
\move(100 60) \rlvec(3 -7) \rlvec(-6 0) \rlvec(3 7)

\move(10 40) \rlvec(10 10) \rlvec(10 -10)
\lpatt(0.3 0.5)
\move(75 40) \rlvec(15 10) \rlvec(15 -10)
\lpatt()
\move(80 40) \rlvec(10 10) \rlvec(10 -10)

\move(20 50) \rlvec(3 -7) \rlvec(-6 0) \rlvec(3 7)
\move(10 40) \rlvec(3 -7) \rlvec(-6 0) \rlvec(3 7)
\move(30 40) \rlvec(3 -7) \rlvec(-6 0) \rlvec(3 7)
\move(90 50) \rlvec(3 -7) \rlvec(-6 0) \rlvec(3 7)
\move(80 40) \rlvec(3 -7) \rlvec(-6 0) \rlvec(3 7)
\move(100 40) \rlvec(3 -7) \rlvec(-6 0) \rlvec(3 7)

\move(55 65) \textref h:R v:C \htext{$\longrightarrow$}
\move(53 68) \textref h:R v:C \htext{\footnotesize{S}}
\move(21.5 65) \textref h:R v:C \htext{\footnotesize{$A$}}
\move(91.5 65) \textref h:R v:C \htext{\footnotesize{$A$}}
\move(31.5 55) \textref h:R v:C \htext{\footnotesize{$C$}}
\move(101.5 55) \textref h:R v:C \htext{\footnotesize{$C$}}

\move(55 45) \textref h:R v:C \htext{$\longrightarrow$}
\move(53 48) \textref h:R v:C \htext{\footnotesize{S}}
\move(21.5 45) \textref h:R v:C \htext{\footnotesize{$A$}}
\move(91.5 45) \textref h:R v:C \htext{\footnotesize{$A$}}
\move(11.5 35) \textref h:R v:C \htext{\footnotesize{$B$}}
\move(81.5 35) \textref h:R v:C \htext{\footnotesize{$B$}}
\move(31.5 35) \textref h:R v:C \htext{\footnotesize{$C$}}
\move(101.5 35) \textref h:R v:C \htext{\footnotesize{$C$}}

} \caption{The shifting operation.}\label{f-sketch-1}
\end{figure}

Figure \ref{f-sketch-1} gives an illustration of the two cases for
the shifting operation, where $A$, $B$ and $C$ denote
subtrees which may be empty, and $S$ stands for the shifting
operation.

The {lifting operation} is meant to transform an even tree $T$ in
position $(i,j)$  to an even tree $T'$ in the position immediately
above $T$. Note that the position $(i,j+1)$ cannot go beyond  the
line $x=2y$.

There are also two cases for the lifting operation. If $T$ contains
a pair of dotted edges, namely, $i=2m+1$, then $T'$ is obtained from
$T$ by moving the pair of edges of the root that are next to the
dotted edges (along with the subtrees attached to these two edges) as outside edges to the last child of the root.

 When $T$ does not contain any dotted edges, namely, $i=2m$,  if $0\leq j\leq m-2$, then $T'$ is
obtained from $T$ by moving the pair of outside edges of the root (along with the subtrees attached to these two edges) as outside edges to the second
child of the root  of $T$, and if $j=m-1$, then let $T'=T$.

\begin{figure}[h,t]
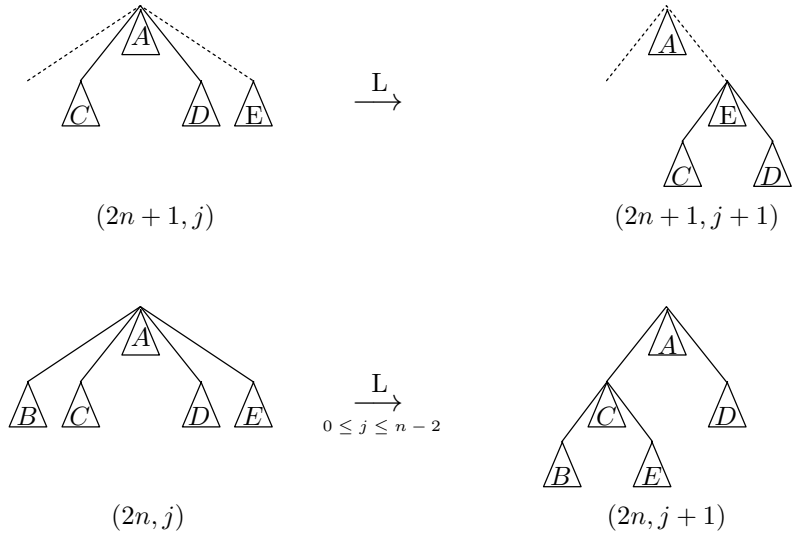

\centertexdraw{
\drawdim mm

\linewd 0.2

\lpatt(0.3 0.5)
\move(5 80) \rlvec(15 10) \rlvec(15 -10)
\move(82 80) \rlvec(8 10) \rlvec(8 -10)
\lpatt()
\move(12 80) \rlvec(8 10) \rlvec(8 -10)
\move(92 72) \rlvec(6 8) \rlvec(6 -8)

\move(20 89.5) \rlvec(2.5 -6) \rlvec(-5 0) \rlvec(2.5 6)
\move(12 80) \rlvec(2.5 -6) \rlvec(-5 0) \rlvec(2.5 6)
\move(28 80) \rlvec(2.5 -6) \rlvec(-5 0) \rlvec(2.5 6)
\move(35 80) \rlvec(2.5 -6) \rlvec(-5 0) \rlvec(2.5 6)

\move(90 89.5) \rlvec(2.5 -6) \rlvec(-5 0) \rlvec(2.5 6)
\move(92 72) \rlvec(2.5 -6) \rlvec(-5 0) \rlvec(2.5 6)
\move(104 72) \rlvec(2.5 -6) \rlvec(-5 0) \rlvec(2.5 6)
\move(98 80) \rlvec(2.5 -6) \rlvec(-5 0) \rlvec(2.5 6)

\move(5 40) \rlvec(15 10) \rlvec(15 -10)
\move(12 40) \rlvec(8 10) \rlvec(8 -10)

\move(82 40) \rlvec(8 10) \rlvec(8 -10)
\move(76 32) \rlvec(6 8) \rlvec(6 -8)

\move(20 49.5) \rlvec(2.5 -6) \rlvec(-5 0) \rlvec(2.5 6)
\move(5 40) \rlvec(2.5 -6) \rlvec(-5 0) \rlvec(2.5 6)
\move(12 40) \rlvec(2.5 -6) \rlvec(-5 0) \rlvec(2.5 6)
\move(28 40) \rlvec(2.5 -6) \rlvec(-5 0) \rlvec(2.5 6)
\move(35 40) \rlvec(2.5 -6) \rlvec(-5 0) \rlvec(2.5 6)

\move(90 49.5) \rlvec(2.5 -6) \rlvec(-5 0) \rlvec(2.5 6)
\move(82 40) \rlvec(2.5 -6) \rlvec(-5 0) \rlvec(2.5 6)
\move(98 40) \rlvec(2.5 -6) \rlvec(-5 0) \rlvec(2.5 6)
\move(76 32) \rlvec(2.5 -6) \rlvec(-5 0) \rlvec(2.5 6)
\move(88 32) \rlvec(2.5 -6) \rlvec(-5 0) \rlvec(2.5 6)

\move(55 77) \textref h:R v:C \htext{$\longrightarrow$}
\move(53 80) \textref h:R v:C \htext{\footnotesize{L}}
\move(21.4 86) \textref h:R v:C \htext{\footnotesize{$A$}}
\move(91.4 85) \textref h:R v:C \htext{\footnotesize{$A$}}
\move(13.4 75.5) \textref h:R v:C \htext{\footnotesize{$C$}}
\move(29.4 75.5) \textref h:R v:C \htext{\footnotesize{$D$}}
\move(36.4 75.5) \textref h:R v:C \htext{\footnotesize{E}}
\move(99.4 75.5) \textref h:R v:C \htext{\footnotesize{E}}
\move(93.4 67.5) \textref h:R v:C \htext{\footnotesize{$C$}}
\move(105.4 67.5) \textref h:R v:C \htext{\footnotesize{$D$}}

\move(55 37) \textref h:R v:C \htext{$\longrightarrow$}
\move(53 40) \textref h:R v:C \htext{\footnotesize{L}}
\move(21.4 46) \textref h:R v:C \htext{\footnotesize{$A$}}
\move(91.4 45) \textref h:R v:C \htext{\footnotesize{$A$}}
\move(6.4 35.5) \textref h:R v:C \htext{\footnotesize{$B$}}
\move(13.4 35.5) \textref h:R v:C \htext{\footnotesize{$C$}}
\move(29.4 35.5) \textref h:R v:C \htext{\footnotesize{$D$}}
\move(36.4 35.5) \textref h:R v:C \htext{\footnotesize{$E$}}
\move(83.4 35.5) \textref h:R v:C \htext{\footnotesize{$C$}}
\move(99.4 35.5) \textref h:R v:C \htext{\footnotesize{$D$}}
\move(77.4 27.5) \textref h:R v:C \htext{\footnotesize{$B$}}
\move(89.4 27.5) \textref h:R v:C \htext{\footnotesize{$E$}}

\move(30 62) \textref h:R v:C \htext{\footnotesize{$(2n+1,j)$}}
\move(105 62) \textref h:R v:C \htext{\footnotesize{$(2n+1,j+1)$}}

\move(26 22) \textref h:R v:C \htext{\footnotesize{$(2n,j)$}}
\move(98 22) \textref h:R v:C \htext{\footnotesize{$(2n,j+1)$}}
\move(60 34) \textref h:R v:C \htext{\tiny{$0\leq j\leq n-2$}}

} \caption{The lifting operation.}\label{f-sketch-2}
\end{figure}

These two cases are illustrated in Figure \ref{f-sketch-2}, where
$A$, $B$, $C$, $D$ and $E$ represent subtrees that may be empty and
$L$ stands for the lifting operation.

Note that in the above process, if $T$ has a pair of dotted edges, then the
degree of the first child of the root must be zero.

Obviously, in the case when $T$  contains a pair of dotted edges,
the lifting operation on $T$ increases the
$r$-index by one, that is, $r(T')=r(T)+1$. It is also clear that the
number of even trees generated at the position $(2n,n)$ equals the
$3$-Catalan numbers. However, it is still necessary to show that the above shifting and
lifting operations do not generate any tree more than once at any
point $(i,j)$.

\begin{figure}[h,t]
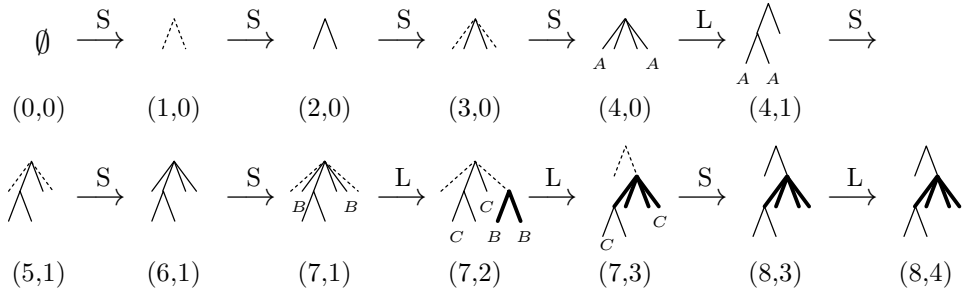

\centertexdraw{
\drawdim mm

\linewd 0.15

\move(15 36) \textref h:R v:C \htext{{$\emptyset$}}
\move(25 36) \textref h:R v:C \htext{$\longrightarrow$}
\move(23 39) \textref h:R v:C \htext{\footnotesize{S}}

\lpatt(0.3 0.5)
\move(30 35) \rlvec(1.5 4) \rlvec(1.5 -4)
\lpatt()
\move(45 36) \textref h:R v:C \htext{$\longrightarrow$}
\move(43 39) \textref h:R v:C \htext{\footnotesize{S}}

\move(50 35) \rlvec(1.5 4) \rlvec(1.5 -4)
\move(65 36) \textref h:R v:C \htext{$\longrightarrow$}
\move(63 39) \textref h:R v:C \htext{\footnotesize{S}}

\lpatt(0.3 0.5)
\move(68.5 35) \rlvec(3 4) \rlvec(3 -4)
\lpatt()
\move(70 35) \rlvec(1.5 4) \rlvec(1.5 -4)
\move(85 36) \textref h:R v:C \htext{$\longrightarrow$}
\move(83 39) \textref h:R v:C \htext{\footnotesize{S}}

\linewd 0.15
\move(88.5 35) \rlvec(3 4) \rlvec(3 -4)
\linewd 0.2
\move(90 35) \rlvec(1.5 4) \rlvec(1.5 -4)
\move(105 36) \textref h:R v:C \htext{$\longrightarrow$}
\move(103 39) \textref h:R v:C \htext{\footnotesize{L}}
\move(96 33) \textref h:R v:C \htext{\tiny{$A$}}
\move(89 33) \textref h:R v:C \htext{\tiny{$A$}}

\move(109 37) \rlvec(1.5 4) \rlvec(1.5 -4)
\move(107.5 33) \rlvec(1.5 4) \rlvec(1.5 -4)
\linewd 0.15
\move(125 36) \textref h:R v:C \htext{$\longrightarrow$}
\move(123 39) \textref h:R v:C \htext{\footnotesize{S}}
\move(108 31) \textref h:R v:C \htext{\tiny{$A$}}
\move(112 31) \textref h:R v:C \htext{\tiny{$A$}}

\lpatt(0.3 0.5)
\move(9.5 16) \rlvec(3 4) \rlvec(3 -4)
\lpatt( )
\move(11 16) \rlvec(1.5 4) \rlvec(1.5 -4)
\move(9.5 12) \rlvec(1.5 4) \rlvec(1.5 -4)
\move(25 15) \textref h:R v:C \htext{$\longrightarrow$}
\move(23 18) \textref h:R v:C \htext{\footnotesize{S}}

\linewd 0.15
\move(28.5 16) \rlvec(3 4) \rlvec(3 -4)
\move(30 16) \rlvec(1.5 4) \rlvec(1.5 -4)
\move(28.5 12) \rlvec(1.5 4) \rlvec(1.5 -4)
\move(45 15) \textref h:R v:C \htext{$\longrightarrow$}
\move(43 18) \textref h:R v:C \htext{\footnotesize{S}}

\linewd 0.15
\lpatt(0.3 0.5)
\move(47 16) \rlvec(4.5 4) \rlvec(4.5 -4)
\lpatt( )
\move(48.5 16) \rlvec(3 4) \rlvec(3 -4)
\move(50 16) \rlvec(1.5 4) \rlvec(1.5 -4)
\move(48.5 12) \rlvec(1.5 4) \rlvec(1.5 -4)
\move(65 15) \textref h:R v:C
\htext{$\longrightarrow$}
\move(63 18) \textref h:R v:C \htext{\footnotesize{L}}
\move(49 14) \textref h:R v:C \htext{\tiny{$B$}}
\move(56 14) \textref h:R v:C \htext{\tiny{$B$}}

\linewd 0.15
\lpatt(0.3 0.5)
\move(67 16) \rlvec(4.5 4) \rlvec(4.5 -4)
\lpatt( )
\move(70 16) \rlvec(1.5 4) \rlvec(1.5 -4)
\move(68.5 12) \rlvec(1.5 4) \rlvec(1.5 -4)
\linewd 0.5
\move(74.5 12) \rlvec(1.5 4) \rlvec(1.5 -4)
\move(85 15) \textref h:R v:C \htext{$\longrightarrow$}
\move(83 18) \textref h:R v:C \htext{\footnotesize{L}}
\move(75 10) \textref h:R v:C \htext{\tiny{$B$}}
\move(79 10) \textref h:R v:C \htext{\tiny{$B$}}
\move(74 14) \textref h:R v:C \htext{\tiny{$C$}}
\move(70 10) \textref h:R v:C \htext{\tiny{$C$}}

\lpatt(0.3 0.5)
\linewd 0.15
\move(90 18) \rlvec(1.5 4) \rlvec(1.5 -4)
\lpatt( )
\move(88.5 10) \rlvec(1.5 4) \rlvec(1.5 -4)
\linewd 0.5
\move(90 14) \rlvec(3 4) \rlvec(3 -4)
\move(91.5 14) \rlvec(1.5 4) \rlvec(1.5 -4)
\move(105 15) \textref h:R v:C \htext{$\longrightarrow$}
\move(103 18) \textref h:R v:C \htext{\footnotesize{S}}
\move(97 12) \textref h:R v:C \htext{\tiny{$C$}}
\move(90 9) \textref h:R v:C \htext{\tiny{$C$}}

\linewd 0.15
\move(110 18) \rlvec(1.5 4) \rlvec(1.5 -4)
\move(108.5 10) \rlvec(1.5 4) \rlvec(1.5 -4)
\linewd 0.5
\move(110 14) \rlvec(3 4) \rlvec(3 -4)
\move(111.5 14) \rlvec(1.5 4) \rlvec(1.5 -4)
\move(125 15) \textref h:R v:C \htext{$\longrightarrow$}
\move(123 18) \textref h:R v:C \htext{\footnotesize{L}}

\linewd 0.15
\move(130 18) \rlvec(1.5 4) \rlvec(1.5 -4)
\move(128.5 10) \rlvec(1.5 4) \rlvec(1.5 -4)
\linewd 0.5
\move(130 14) \rlvec(3 4) \rlvec(3 -4)
\move(131.5 14) \rlvec(1.5 4) \rlvec(1.5 -4)

\move(17 27) \textref h:R v:C \htext{\footnotesize{(0,0)}}
\move(35 27) \textref h:R v:C \htext{\footnotesize{(1,0)}}
\move(55 27) \textref h:R v:C \htext{\footnotesize{(2,0)}}
\move(75 27) \textref h:R v:C \htext{\footnotesize{(3,0)}}
\move(95 27) \textref h:R v:C \htext{\footnotesize{(4,0)}}
\move(115 27) \textref h:R v:C \htext{\footnotesize{(4,1)}}

\move(17 5) \textref h:R v:C \htext{\footnotesize{(5,1)}}
\move(35 5) \textref h:R v:C \htext{\footnotesize{(6,1)}}
\move(55 5) \textref h:R v:C \htext{\footnotesize{(7,1)}}
\move(75 5) \textref h:R v:C \htext{\footnotesize{(7,2)}}
\move(95 5) \textref h:R v:C \htext{\footnotesize{(7,3)}}
\move(115 5) \textref h:R v:C \htext{\footnotesize{(8,3)}}
\move(135 5) \textref h:R v:C \htext{\footnotesize{(8,4)}}

} \caption{The actions of the shifting and lifting
operations.}\label{f-eventree-2}
\end{figure}

We use induction on $i$ and $j$ to complete the proof. It is easy to check the
statement  holds for small values of $(i,j)$, see
Figure \ref{f-eventree-2}. Assume that the claim holds for the positions
$(i-1, j)$ and $(i,j-1)$.

Suppose that $T_1'$ and $T_2'$ are two even
trees in the position $(i,j)$ that are obtained from $T_1$ and
$T_2$, respectively. There are three cases for the positions of $T_1$
and $T_2$.

Case 1. Both $T_1$ and $T_2$ are in the position $(i-1,j)$. In this
case,  we have $i\geq 1$. Clearly, $T_1'$ and $T_2'$ are obtained from
$T_1$ and $T_2$  via the shifting operations. It is also easy
to see  that
$T_1'$ and $T_2'$ are distinct.

Case 2. Both $T_1$ and $T_2$ are in the position $(i,j-1)$, where
$j\geq 1$. In this case, $T_1'$ and $T_2'$ are obtained from $T_1$
and $T_2$ via the lifting operations. As can be seen from Figure
\ref{f-sketch-2}, if $T_1$ and $T_2$ has dotted edges, then we
consider the  subtrees $A_1, C_1, D_1, E_1$ and $A_2, C_2, D_2,
E_2$ of $T_1$ and $T_2$. Bear in mind that these subtrees may be
empty. Since $T_1$ and $T_2$ are distinct, the corresponding
subtrees cannot be all identical. After the lifting operation, the four
subtrees are moved to different positions in $T_1'$ and $T_2'$. It
is easily seen that $T_1'$ and $T_2'$ are distinct.

If $T_1$ and $T_2$ do  not have any dotted edges, we may consider
the five subtrees $A, B, C, D, E$ as shown in Figure
\ref{f-sketch-2}. The same argument yields that $T_1'$ and $T_2'$
are distinct.

Case 3.
   $T_1$ is in the position $(i-1,j)$ and $T_2$ is in the
position $(i,j-1)$. In this case, we have $i,j\geq 1$. Then $T_1'$
is obtained from $T_1$ via a shifting operation and $T_2'$ is
obtained from $T_2$ via a lifting operation.

If $T_1$ has dotted edges, then the first child of the root is a
leaf and the shifting operation transforms the dotted edges to
regular edges. Therefore, the first child of root of $T_1'$ is also
a leaf. On the other hand, $T_2$ does not have any dotted edges. Applying the lifting operation to $T_2$, the first child of the root
of $T_2'$ becomes an internal vertex. Hence  $T_1'$ and $T_2'$ are distinct.

If $T_1$ does not have any dotted edges, the shifting operation
adds two dotted edges as outside edges to the root.
Moreover, the last child of the root of $T_1'$  is a leaf. Meanwhile,
since $T_2$ has dotted edges, applying the lifting operation to
$T_2$,  the last child of the root of $T_2'$ becomes an internal vertex. Thus, $T_1'$ and $T_2'$ are distinct.

Since the number of even trees with $2n$ edges
equals the $3$-Catalan number $T_n$, the above construction produces all even trees with $2n$ edges at the position $(2n,n)$. Hence we obtain a one-to-one correspondence.
This completes the proof.  \qed

\begin{figure}[h,t]
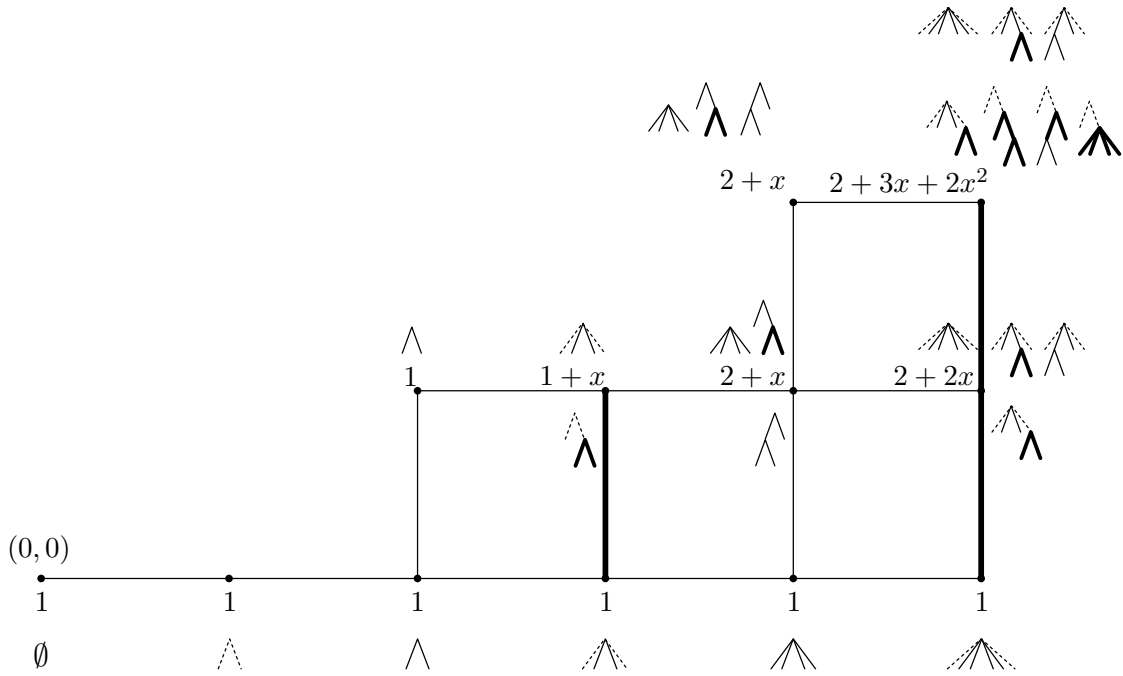

\centertexdraw{
\drawdim mm

\linewd 0.2
\move(10 15) \rlvec(125 0) \rlvec(0 50)
\move(60 15) \rlvec(0 25) \rlvec(75 0)
\move(110 15) \rlvec(0 50) \rlvec(25 0)
\linewd 0.8
\move(85 15) \rlvec(0 25)
\move(135 15) \rlvec(0 50)

\move(10 15) \fcir f:0 r:0.5
\move(35 15) \fcir f:0 r:0.5
\move(60 15) \fcir f:0 r:0.5
\move(85 15) \fcir f:0 r:0.5
\move(110 15) \fcir f:0 r:0.5
\move(135 15) \fcir f:0 r:0.5

\move(60 40) \fcir f:0 r:0.5
\move(85 40) \fcir f:0 r:0.5
\move(110 40) \fcir f:0 r:0.5
\move(135 40) \fcir f:0 r:0.5

\move(110 65) \fcir f:0 r:0.5
\move(135 65) \fcir f:0 r:0.5

\linewd 0.15
\move(11 5) \textref h:R v:C \htext{{$\emptyset$}}

\lpatt(0.3 0.5)
\move(33.5 3) \rlvec(1.5 4) \rlvec(1.5 -4)

\lpatt()
\move(58.5 3) \rlvec(1.5 4) \rlvec(1.5 -4)

\lpatt(0.3 0.5)
\move(82 3) \rlvec(3 4) \rlvec(3 -4)
\lpatt()
\move(83.5 3) \rlvec(1.5 4) \rlvec(1.5 -4)

\move(107 3) \rlvec(3 4) \rlvec(3 -4)
\move(108.5 3) \rlvec(1.5 4) \rlvec(1.5 -4)

\lpatt(0.3 0.5)
\move(130.5 3) \rlvec(4.5 4) \rlvec(4.5 -4)
\lpatt()
\move(132 3) \rlvec(3 4) \rlvec(3 -4)
\move(133.5 3) \rlvec(1.5 4) \rlvec(1.5 -4)

\move(58 45) \rlvec(1.3 3.5) \rlvec(1.3 -3.5)

\lpatt(0.3 0.5)
\move(79 45) \rlvec(3 4) \rlvec(3 -4)
\lpatt()
\move(80.5 45) \rlvec(1.5 4) \rlvec(1.5 -4)

\lpatt(0.3 0.5)
\move(79.7 33.5) \rlvec(1.3 3.5) \rlvec(1.3 -3.5)
\lpatt( )
\linewd 0.5
\move(81 30) \rlvec(1.3 3.5) \rlvec(1.3 -3.5)

\linewd 0.15
\move(99 45) \rlvec(2.6 3.5) \rlvec(2.6 -3.5)
\move(100.3 45) \rlvec(1.3 3.5) \rlvec(1.3 -3.5)

\move(104.7 48.5) \rlvec(1.3 3.5) \rlvec(1.3 -3.5)
\linewd 0.5
\move(106 45) \rlvec(1.3 3.5) \rlvec(1.3 -3.5)

\linewd 0.15
\lpatt()
\move(106.3 33.5) \rlvec(1.3 3.5) \rlvec(1.3 -3.5)
\move(105 30) \rlvec(1.3 3.5) \rlvec(1.3 -3.5)

\linewd 0.15
\lpatt(0.3 0.5)
\move(126.7 45.5) \rlvec(3.9 3.5) \rlvec(3.9 -3.5)
\lpatt( )
\move(128 45.5) \rlvec(2.6 3.5) \rlvec(2.6 -3.5)
\move(129.3 45.5) \rlvec(1.3 3.5) \rlvec(1.3 -3.5)

\lpatt(0.3 0.5)
\move(136.4 45.5) \rlvec(2.6 3.5) \rlvec(2.6 -3.5)
\lpatt( )
\move(137.7 45.5) \rlvec(1.3 3.5) \rlvec(1.3 -3.5)
\linewd 0.5
\move(139 42) \rlvec(1.3 3.5) \rlvec(1.3 -3.5)

\linewd 0.15
\lpatt(0.3 0.5)
\move(143.4 45.5) \rlvec(2.6 3.5) \rlvec(2.6 -3.5)
\lpatt( )
\move(144.7 45.5) \rlvec(1.3 3.5) \rlvec(1.3 -3.5)
\move(143.4 42) \rlvec(1.3 3.5) \rlvec(1.3 -3.5)

\lpatt(0.3 0.5)
\move(136.4 34.5) \rlvec(2.6 3.5) \rlvec(2.6 -3.5)
\lpatt( )
\move(137.7 34.5) \rlvec(1.3 3.5) \rlvec(1.3 -3.5)
\linewd 0.5
\move(140.3 31) \rlvec(1.3 3.5)  \rlvec(1.3 -3.5)

\linewd 0.15
\lpatt()
\move(90.8 74.5) \rlvec(2.6 3.5) \rlvec(2.6 -3.5)
\move(92.1 74.5) \rlvec(1.3 3.5) \rlvec(1.3 -3.5)

\lpatt()
\move(97.1 77.5) \rlvec(1.3 3.5) \rlvec(1.3 -3.5)
\linewd 0.5
\move(98.4 74) \rlvec(1.3 3.5) \rlvec(1.3 -3.5)

\linewd 0.15
\lpatt()
\move(104.3 77.5) \rlvec(1.3 3.5) \rlvec(1.3 -3.5)
\move(103 74) \rlvec(1.3 3.5) \rlvec(1.3 -3.5)

\linewd 0.15
\lpatt(0.3 0.5)
\move(126.7 87.5) \rlvec(3.9 3.5) \rlvec(3.9 -3.5)
\lpatt( )
\move(128 87.5) \rlvec(2.6 3.5) \rlvec(2.6 -3.5)
\move(129.3 87.5) \rlvec(1.3 3.5) \rlvec(1.3 -3.5)

\lpatt(0.3 0.5)
\move(136.4 87.5) \rlvec(2.6 3.5) \rlvec(2.6 -3.5)
\lpatt( )
\move(137.7 87.5) \rlvec(1.3 3.5) \rlvec(1.3 -3.5)
\linewd 0.5
\move(139 84) \rlvec(1.3 3.5) \rlvec(1.3 -3.5)

\linewd 0.15
\lpatt(0.3 0.5)
\move(143.4 87.5) \rlvec(2.6 3.5) \rlvec(2.6 -3.5)
\lpatt( )
\move(144.7 87.5) \rlvec(1.3 3.5) \rlvec(1.3 -3.5)
\move(143.4 84) \rlvec(1.3 3.5) \rlvec(1.3 -3.5)

\linewd 0.15
\lpatt(0.3 0.5)
\move(127.8 75) \rlvec(2.6 3.5) \rlvec(2.6 -3.5)
\lpatt( )
\move(129.1 75) \rlvec(1.3 3.5) \rlvec(1.3 -3.5)
\linewd 0.5
\move(131.7 71.5) \rlvec(1.3 3.5) \rlvec(1.3 -3.5)

\linewd 0.15
\lpatt(0.3 0.5)
\move(135.4 77) \rlvec(1.3 3.5) \rlvec(1.3 -3.5)
\lpatt( )
\linewd 0.5
\move(136.7 73.5) \rlvec(1.3 3.5) \rlvec(1.3 -3.5)
\linewd 0.5
\move(138.1 70) \rlvec(1.3 3.5) \rlvec(1.3 -3.5)

\linewd 0.15
\lpatt(0.3 0.5)
\move(142.4 77) \rlvec(1.3 3.5) \rlvec(1.3 -3.5)
\lpatt( )
\linewd 0.5
\move(143.7 73.5) \rlvec(1.3 3.5)  \rlvec(1.3 -3.5)
\linewd 0.15
\move(142.4 70) \rlvec(1.3 3.5) \rlvec(1.3 -3.5)

\linewd 0.15
\lpatt(0.3 0.5)
\move(148.1 75) \rlvec(1.3 3.5) \rlvec(1.3 -3.5)
\linewd 0.5
\lpatt( )
\move(148.1 71.5) \rlvec(2.6 3.5)\rlvec(2.6 -3.5)
\move(149.4 71.5) \rlvec(1.3 3.5) \rlvec(1.3 -3.5)

\move(14 19) \textref h:R v:C \htext{\small{$(0,0)$}}

\move(11 12) \textref h:R v:C \htext{\small{$1$}}
\move(36 12) \textref h:R v:C \htext{\small{$1$}}
\move(61 12) \textref h:R v:C \htext{\small{$1$}}
\move(86 12) \textref h:R v:C \htext{\small{$1$}}
\move(111 12) \textref h:R v:C \htext{\small{$1$}}
\move(136 12) \textref h:R v:C \htext{\small{$1$}}

\move(60 42) \textref h:R v:C \htext{\small{$1$}}
\move(85 42) \textref h:R v:C \htext{\small{$1+x$}}
\move(109 42) \textref h:R v:C \htext{\small{$2+x$}}
\move(134 42) \textref h:R v:C \htext{\small{$2+2x$}}

\move(109 68) \textref h:R v:C \htext{\small{$2+x$}}
\move(136 68) \textref h:R v:C \htext{\small{$2+3x+2x^2$}}

} \caption{The generation of the polynomials $R_n(x)$.}
\label{f-eventree-1}
\end{figure}

For example, Figure \ref{f-eventree-2} gives a sequence of shifting
and lifting operations initially acting on the empty even tree in
the origin $(0,0)$. Figure \ref{f-eventree-1} gives the first
few steps to generate even trees.

\vskip 5pt
 \noindent \textbf{Acknowledgments.}  This work was
supported by the 973 Project, the PCSIRT Project of the Ministry of
Education, and the National Science Foundation of China.

\end{document}